\documentclass[10pt]{amsart}

\usepackage{amsmath,amsfonts,amssymb,amsthm,amscd,multicol,epsfig,psfrag}

\def\C{\mathbb{C}}
\def\Z{\mathbb{Z}}

\def\Q{\mathbb{Q}}

\def\id{\mathbf{1}}

\def\g{\ensuremath{\mathfrak{g}}}

\def\V{\mathbf{V}}

\def\S{\mathbf{S}}

\def\P{\mathcal{P}}

\newcommand{\comment}[1]{}

\DeclareMathOperator{\Hom}{Hom}
\DeclareMathOperator{\rank}{rank}
\DeclareMathOperator{\End}{End}
\DeclareMathOperator{\Aut}{Aut}

\DeclareMathOperator{\tr}{tr}
\DeclareMathOperator{\Irr}{Irr}

\newtheorem{theo}{Theorem}[section]
\newtheorem{prop}[theo]{Proposition}
\newtheorem{lem}[theo]{Lemma}

\newtheorem{defin}[theo]{Definition}
\newtheorem{example}[theo]{Example}
\newtheorem*{rem*}{Remark}

\numberwithin{equation}{section}

\begin{document}
\title{Finite-dimensional algebras and quivers}
\author{Alistair Savage}
\address{Fields Institute and University of Toronto \\
Toronto, Ontario \\ Canada} \email{alistair.savage@aya.yale.edu}
\thanks{This research was supported in part by the Natural
Sciences and Engineering Research Council (NSERC) of Canada}
\date{May 4, 2005}

\begin{abstract}
This is an overview article on finite-dimensional algebras and
quivers, written for the {\em Encyclopedia of Mathematical Physics}.
We cover path algebras, Ringel-Hall algebras and the quiver
varieties of Lusztig and Nakajima.
\end{abstract}

\maketitle

\begin{multicols}{2}

\section{Introduction}
Algebras and their representations are ubiquitous in mathematics.
It turns out that representations of finite-dimensional algebras
are intimately related to quivers, which are simply oriented
graphs.  Quivers arise naturally in many areas of mathematics,
including representation theory, algebraic and differential
geometry, Kac-Moody algebras, and quantum groups. In this article,
we give a brief overview of some of these topics. We start by
giving the basic definitions of associative algebras and their
representations.  We then introduce quivers and their
representation theory, mentioning the connection to the
representation theory of associative algebras.  We also discuss in
some detail the relationship between quivers and the theory of Lie
algebras.


\section{Associative algebras}

An \emph{algebra} is a vector space $A$ over a field $k$ equipped
with a multiplication which is distributive and such that
\[
a(xy) = (ax)y = x(ay),\ \forall\ a \in k,\ x,y \in A.
\]
When we wish to make the field explicit, we call $A$ a
$k$-algebra.  An algebra is \emph{associative} if $(xy)z = x(yz)$
for all $x,y,z \in A$. $A$ has a \emph{unit}, or
\emph{multiplicative identity}, if it contains an element $\id_A$
such that $\id_A x = x \id_A = x$ for all $x \in A$. From now on,
we will assume all algebras are associative with unit.  $A$ is
said to be \emph{commutative} if $xy=yx$ for all $x,y \in A$ and
\emph{finite-dimensional} if the underlying vector space of $A$ is
finite-dimensional.

A vector subspace $I$ of $A$ is called a \emph{left (resp. right)
ideal} if $xy \in I$ for all $x \in A$, $y \in I$ (resp. $x \in
I$, $y \in A$).  If $I$ is both a right and a left ideal, it is
called a \emph{two-sided} ideal of $A$.  If $I$ is a two-sided
ideal of $A$, then the factor space $A/I$ is again an algebra.

An \emph{algebra homomorphism} is a linear map $f : A_1 \to A_2$
between two algebras such that
\begin{gather*}
f(\id_{A_1}) = \id_{A_2},\ \text{and} \\
f(xy) = f(x) f(y),\ \forall\ x,y \in A.
\end{gather*}
A \emph{representation} of an algebra $A$ is an algebra
homomorphism $\rho : A \to \End_k (V)$ for a $k$-vector space $V$.
Here $\End_k (V)$ is the space of endomorphisms of the vector
space $V$ with multiplication given by composition.  Given a
representation of an algebra $A$ on a vector space $V$, we may
view $V$ as an $A$-module with the action of $A$ on $V$ given by
\[
a \cdot v = \rho(a)v,\ a \in A, v \in V.
\]

A \emph{morphism} $\psi : V \to W$ of two $A$-modules (or
equivalently, representations of $A$), is a linear map commuting
with the action of $A$.  That is, it is a linear map satisfying
\[
a \cdot \psi(v) = \psi(a \cdot v),\ \forall\ a \in A,\, v \in V.
\]

Let $G$ be a commutative monoid (a set with an associative
multiplication and a unit element).  A \emph{$G$-graded
$k$-algebra} is a $k$-algebra which can be expressed as a direct
sum $A=\oplus_{g \in G} A_g$ such that $aA_g \subset A_g$ for all
$a \in k$ and $A_{g_1} A_{g_2} \subset A_{g_1 + g_2}$ for all
$g_1, g_2 \in G$. A morphism $\psi : A \to B$ of $G$-graded
algebras is a $k$-algebra morphism respecting the grading, that
is, satisfying $\psi(A_g) \subset B_g$ for all $g \in G$.


\section{Quivers and path algebras}

A \emph{quiver} is simply an oriented graph.  More precisely, a
quiver is a pair $Q=(Q_0,Q_1)$ where $Q_0$ is a finite set of
vertices and $Q_1$ is a finite set of arrows (oriented edges)
between them.  For $a \in Q_1$, we let $h(a)$ denote the
\emph{head} of $a$ and $t(a)$ denote the \emph{tail} of $a$.  A
\emph{path} in $Q$ is a sequence $x = \rho_1 \rho_2 \dots \rho_m$
of arrows such that $h(\rho_{i+1}) = t(\rho_i)$ for $1 \le i \le
m-1$.  We let $t(x) = t(\rho_m)$ and $h(\rho) = h(\rho_1)$ denote
the initial and final vertices of the path $x$.  For each vertex
$i \in Q_0$, we let $e_i$ denote the trivial path which starts and
ends at the vertex $i$.

Fix a field $k$.  The \emph{path algebra} $kQ$ associated to a
quiver $Q$ is the $k$-algebra whose underlying vector space has
basis the set of paths in $Q$, and with the product of paths given
by concatenation.  Thus, if $x = \rho_1 \dots \rho_m$ and $y =
\sigma_1 \dots \sigma_n$ are two paths, then $xy = \rho_1 \dots
\rho_m \sigma_1 \dots \sigma_n$ if $h(y) = t(x)$ and $xy=0$
otherwise.  We also have
\begin{gather*}
e_i e_j = \begin{cases} e_i \text{ if } i = j \\
0 \text{ if } i \ne j \end{cases},
\end{gather*}
\begin{gather*}
e_i x = \begin{cases} x \text{ if } h(x) = i \\
0 \text{ if } h(x) \ne i \end{cases},
\end{gather*}
\begin{gather*}
x e_i = \begin{cases} x \text{ if } t(x) = i \\
0 \text{ if } t(x) \ne i \end{cases},
\end{gather*}
for $x \in kQ$. This multiplication is associative.  Note that
$e_i A$ and $A e_i$ have bases given by the set of paths ending
and starting at $i$ respectively.  The path algebra has a unit
given by $\sum_{i \in Q_0} e_i$.

\begin{example}
\label{example-3quiver}
Let $Q$ be the following quiver.
\[
\psfrag{r}{$\rho$} \psfrag{s}{$\sigma$} \psfrag{l}{$\lambda$}
\includegraphics[width=0.3\textwidth]{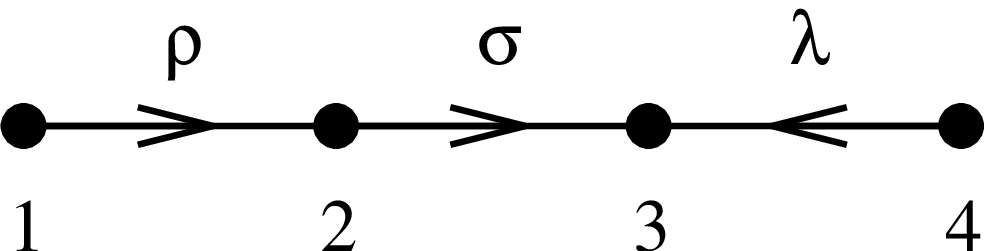}
\]
Then $kQ$ has a basis given by the set of paths $\{e_1,e_2,e_3,
e_4, \rho, \sigma, \lambda, \sigma \rho\}$.  Some sample products
are $\rho \sigma = 0$, $\lambda \lambda = 0$, $\lambda \sigma =
0$, $e_3 \sigma = \sigma e_2 = \sigma$, $e_2 \sigma = 0$.
\end{example}

\begin{example}
Let $Q$ be the following quiver (the so-called \emph{Jordan
quiver}).
\[
\psfrag{r}{$\rho$}
\includegraphics[height=0.15\textwidth]{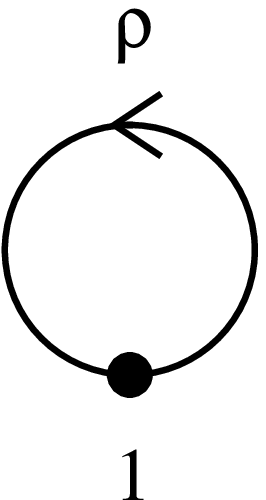}
\]
Then $kQ \cong k[t]$, the algebra of polynomials in one variable.
\end{example}

Note that the path algebra $kQ$ is finite-dimensional if and only if
$Q$ has no oriented cycles (paths with the same head and tail
vertex).

\begin{example}
Let $Q$ be the following quiver.
\\
\[
\psfrag{n2}{$n\! -\! 2$} \psfrag{n1}{$n\! -\! 1$} \psfrag{n}{$n$}
\includegraphics[width=0.4\textwidth]{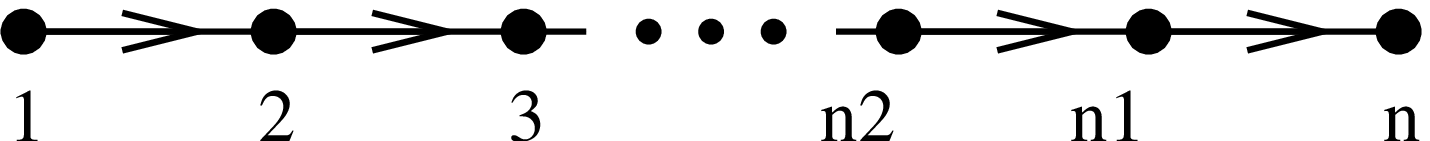}
\]
Then for every $1 \le i \le j \le n$, there is a unique path from
$i$ to $j$.  Let $f : kQ \to M_n(k)$ be the linear map from the
path algebra to the $n \times n$ matrices with entries in the
field $k$ that sends the unique path from $i$ to $j$ to the matrix
$E_{ji}$ with $(j,i)$ entry 1 and all other entries zero. Then one
can show that $f$ is an isomorphism onto the algebra of lower
triangular matrices.
\end{example}

\section{Representations of quivers}

Fix a field $k$.  A representation of a quiver $Q$ is an
assignment of a vector space to each vertex and to each arrow a
linear map between the vector spaces assigned to its tail and
head.  More precisely, a representation $V$ of $Q$ is a collection
\[
\{V_i\ |\ i \in Q_0\}
\]
of finite-dimensional $k$-vector spaces together with a collection
\[
\{V_\rho\, :\, V_{t(\rho)} \to V_{h(\rho)}\ |\ \rho \in Q_1\}
\]
of $k$-linear maps.  Note that a representation $V$ of a quiver
$Q$ is equivalent to a representation of the path algebra $kQ$.
The \emph{dimension} of $V$ is the map $d_V\, :\, Q_0 \to \Z_{\ge
0}$ given by $d_V(i) = \dim V_i$ for $i \in Q_0$.

If $V$ and $W$ are two representations of a quiver $Q$, then a
\emph{morphism} $\psi : V \to W$ is a collection of $k$-linear
maps
\[
\{\psi_i\, :\, V_i \to W_i\ |\ i \in Q_0\}
\]
such that
\[
W_\rho \psi_{t(\rho)} = \psi_{h(\rho)} V_\rho,\ \forall\ \rho \in
Q_1.
\]

\begin{prop}
Let $A$ be a finite-dimensional $k$-algebra.  Then the category of
representations of $A$ is equivalent to the category of
representations of the algebra $kQ/I$ for some quiver $Q$ and some
two-sided ideal $I$ of $kQ$.
\end{prop}
It is for this reason that the study of finite-dimensional
associative algebras is intimately related to the study of
quivers.

We define the \emph{direct sum} $V \oplus W$ of two
representations $V$ and $W$ of a quiver $Q$ by
\[
(V \oplus W)_i = V_i \oplus W_i,\ i \in Q_0
\]
and $(V \oplus W)_\rho : V_{t(\rho)} \oplus W_{t(\rho)} \to
V_{h(\rho)} \oplus W_{h(\rho)}$ by
\[
(V \oplus W)_\rho ((v,w)) = (V_\rho(v),W_\rho(w))
\]
for $v \in V_{t(\rho)}$, $w \in W_{t(\rho)}$, $\rho \in Q_1$.  A
representation $V$ is \emph{trivial} if $V_i=0$ for all $i \in Q_0$
and \emph{simple} if its only subrepresentations are the zero
representation and $V$ itself. We say that $V$ is
\emph{decomposable} if it is isomorphic to $W \oplus U$ for some
nontrivial representations $W$ and $U$. Otherwise, we call $V$
\emph{indecomposable}.  Every representation of a quiver has a
decomposition into indecomposable representations that is unique up
to isomorphism and permutation of the components. Thus, to classify
all representations of a quiver, it suffices to classify the
indecomposable representations.

\begin{example}
Let $Q$ be the following quiver.
\[
\psfrag{r}{$\rho$}
\includegraphics[width=0.1\textwidth]{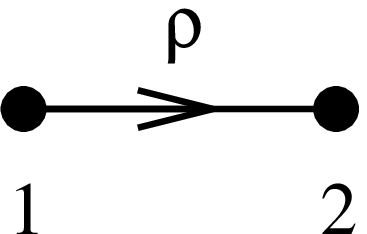}
\]
Then $Q$ has three indecomposable representations $U$, $V$ and $W$
given by:
\begin{gather*}
U_1 = k,\ U_2 = 0,\ U_\rho = 0, \\
V_1 = 0,\ V_2 = k,\ V_\rho = 0, \\
W_1 = k,\ W_2 = k,\ W_\rho = 1.
\end{gather*}
Then any representation $Z$ of $Q$ is isomorphic to
\[
Z \cong U^{d_1 - r} \oplus V^{d_2 - r} \oplus W^r,
\]
where $d_1 = \dim Z_1$, $d_2 = \dim Z_2$, $r = \rank Z_\rho$.
\end{example}

\begin{example}
Let $Q$ be the Jordan quiver.  Then representations $V$ of $Q$ are
classified up to isomorphism by the Jordan normal form of $V_\rho$
where $\rho$ is the single arrow of the quiver.  Indecomposable
representations correspond to single Jordan blocks.  These are
parametrized by a discrete parameter $n$ (the size of the block)
and a continuous parameter $\lambda$ (the eigenvalue of the
block).
\end{example}

A quiver is said to be of \emph{finite type} if it has only
finitely many indecomposable representations (up to isomorphism).
If a quiver has infinitely many isomorphism classes but they can
be split into families, each parametrized by a single continuous
parameter, then we say the quiver is of \emph{tame} (or
\emph{affine}) \emph{type}. If a quiver is of neither finite nor
tame type, it is of \emph{wild type}.  It turns out that there is
a rather remarkable relationship between the classification of
quivers and their representations and the theory of Kac-Moody
algebras.

The \emph{Euler form} or \emph{Ringel form} of a quiver $Q$ is
defined to be the asymmetric bilinear form on $\Z^{Q_0}$ given by
\[
\left<\alpha,\beta\right> = \sum_{i \in Q_0} \alpha(i) \beta(i) -
\sum_{\rho \in Q_1} \alpha(t(\rho)) \beta(h(\rho)).
\]
In the standard coordinate basis of $\Z^{Q_0}$, the Euler form is
represented by the matrix $E=(a_{ij})$ where
\[
a_{ij} = \delta_{ij} - \#\{\rho \in Q_1\ |\ t(\rho) = i,\, h(\rho)
= j\}.
\]
Here $\delta_{ij}$ is the Kronecker delta symbol.  We define the
\emph{Cartan form} of the quiver $Q$ to be the symmetric bilinear
form given by
\[
(\alpha,\beta) = \left<\alpha, \beta\right> + \left<\beta, \alpha
\right>.
\]
Note that the Cartan form is independent of the orientation of the
arrows in $Q$.  In the standard coordinate basis of $\Z^{Q_0}$,
the Cartan form is represented by the \emph{Cartan matrix}
$C=(c_{ij})$ where $c_{ij} = a_{ij} + a_{ji}$.

\begin{example}
For the quiver in Example~\ref{example-3quiver}, the Euler matrix
is
\[
E = \begin{pmatrix} 1 & -1 & 0 & 0 \\
0 & 1 & -1 & 0 \\
0 & 0 & 1 & 0 \\
0 & 0 & -1 & 1 \end{pmatrix}
\]
and the Cartan matrix is
\[
C = \begin{pmatrix} 2 & -1 & 0 & 0 \\
-1 & 2 & -1 & 0 \\
0 & -1 & 2 & -1 \\
0 & 0 & -1 & 2 \end{pmatrix}.
\]
\end{example}

The \emph{Tits form} $q$ of a quiver $Q$ is defined by
\[
q(\alpha) = \left<\alpha, \alpha\right> = \frac{1}{2}(\alpha,
\alpha).
\]
It is known that the number of continuous parameters describing
representations of dimension $\alpha$ for $\alpha \ne 0$ is
greater than or equal to $1 - q(\alpha)$.

Let $\g$ be the Kac-Moody algebra associated to the Cartan matrix
of a quiver $Q$.  By forgetting the orientation of the arrows of
$Q$, we obtain the underlying (undirected) graph.  This is the
Dynkin graph of $\g$.  Associated to $\g$ is a root system and a
set of simple roots $\{\alpha_i\ |\ i \in Q_0\}$ indexed by the
vertices of the Dynkin graph.

\begin{theo}[Gabriel's Theorem]
\
\begin{enumerate}
\item A quiver is of finite type it and only if the underlying
graph is a union of Dynkin graphs of type $A$, $D$, or $E$.

\item A quiver is of tame type if and only if the underlying graph
is a union of Dynkin graphs of type $A$, $D$, or $E$ and extended
Dynkin graphs of type $\widehat A$, $\widehat D$, or $\widehat E$
(with at least one extended Dynkin graph).

\item The isomorphism classes of indecomposable representations of
a quiver $Q$ of finite type are in one-to-one correspondence with
the positive roots of the root system associated to the underlying
graph of $Q$.  The correspondence is given by
\begin{equation*}
V \mapsto \sum_{i \in Q_0} d_V(i) \alpha_i.
\end{equation*}
\end{enumerate}
\end{theo}

The Dynkin graphs of type $A$, $D$, and $E$ are as follows.
\\
\begin{center}
\psfrag{A}{$A_n$} \psfrag{D}{$D_n$} \psfrag{E6}{$E_6$}
\psfrag{E7}{$E_7$} \psfrag{E8}{$E_8$}
\includegraphics[width=0.4\textwidth]{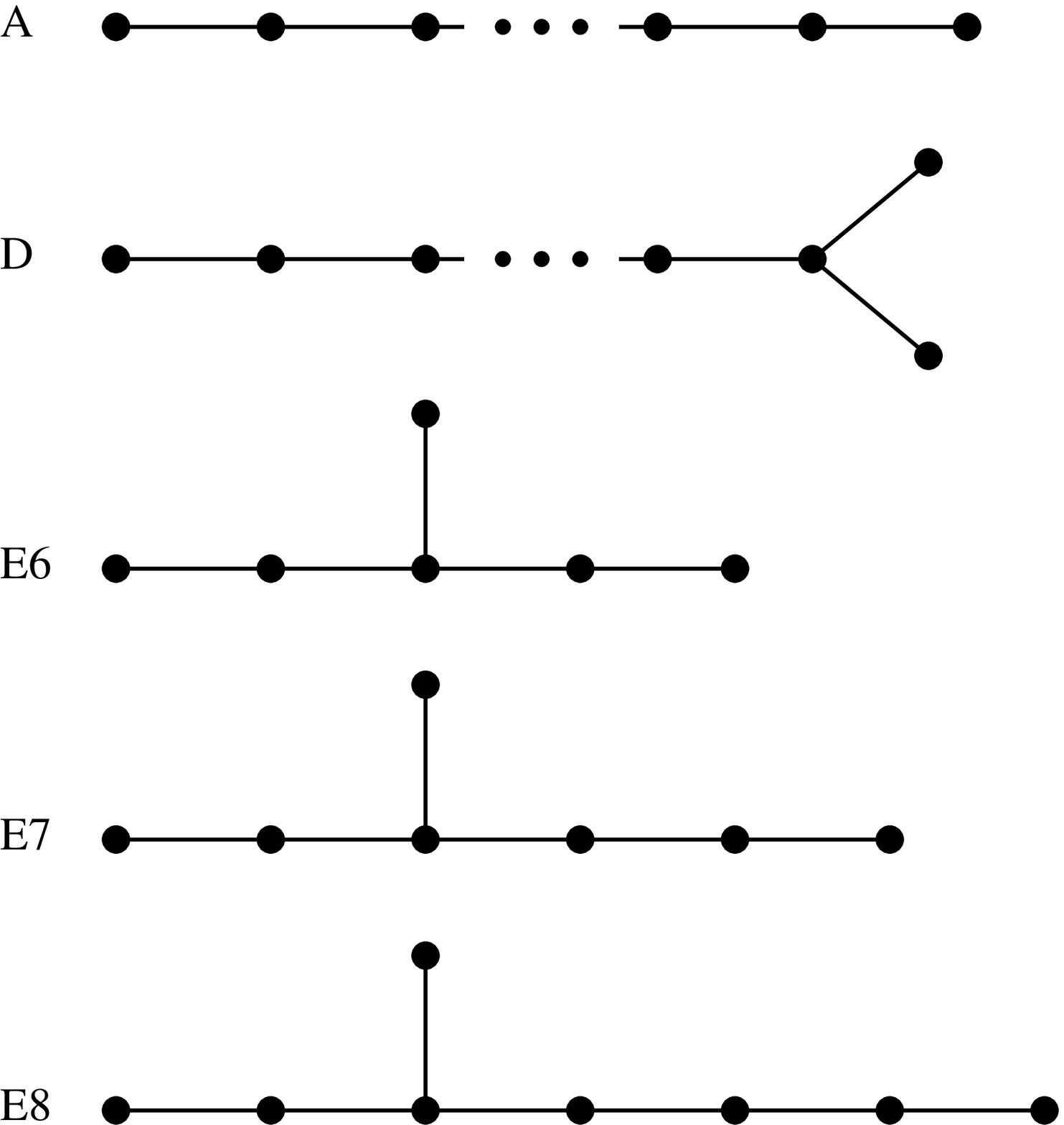}
\end{center}
\ \\
Here the subscript indicates the number of vertices in the
graph.

The extended Dynkin graphs of type $\widehat A$, $\widehat D$, and
$\widehat E$ are as follows.
\begin{center}
\psfrag{Ah}{${\widehat A}_n$} \psfrag{Dh}{${\widehat D}_n$}
\psfrag{E6h}{${\widehat E}_6$} \psfrag{E7h}{${\widehat E}_7$}
\psfrag{E8h}{${\widehat E}_8$}
\includegraphics[width=0.4\textwidth]{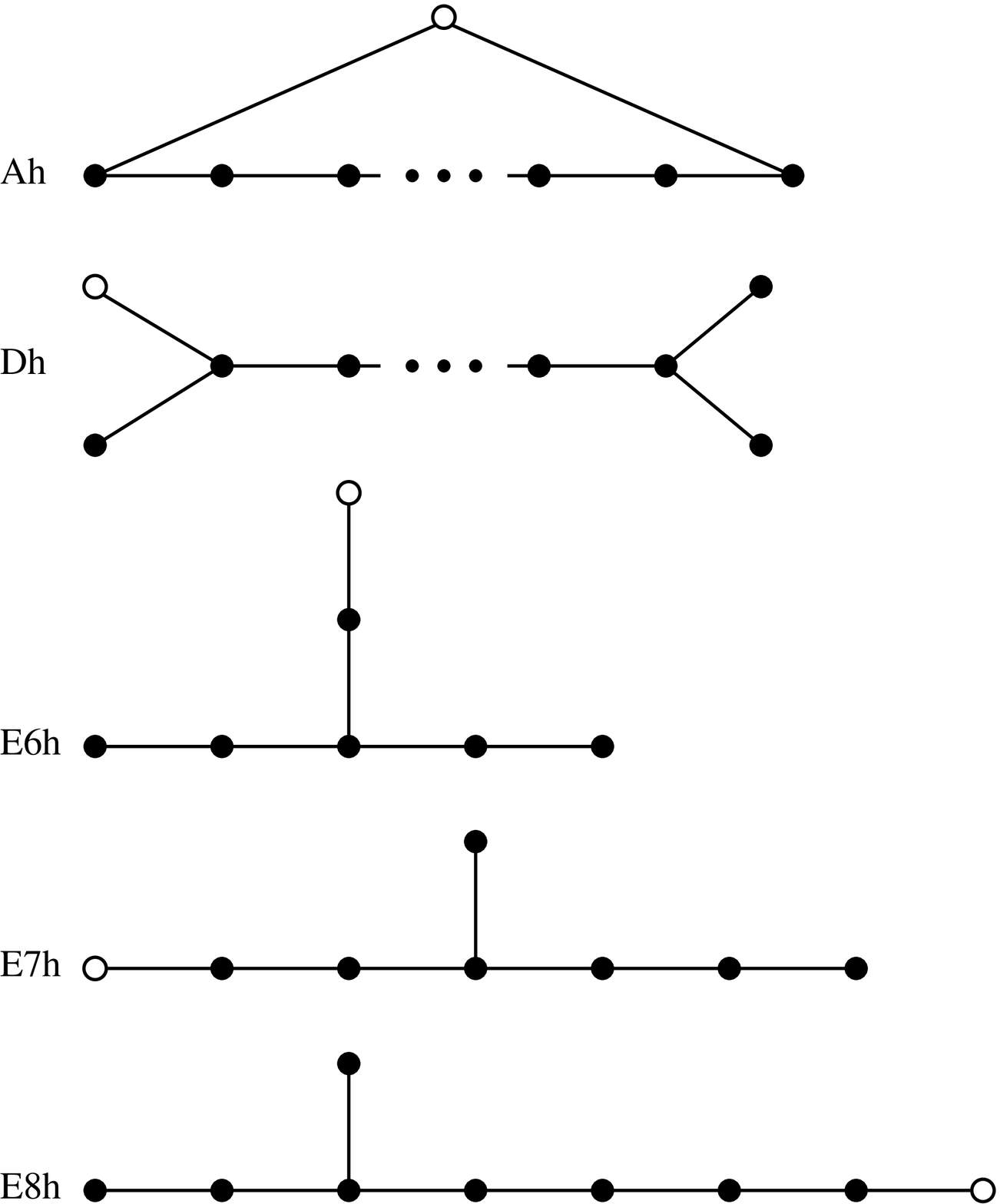}
\end{center}
\ \\
Here we have used an open dot to denote the vertex that was added
to the corresponding Dynkin graph of type $A$, $D$, or $E$.

\begin{theo}[Kac's Theorem]
Let $Q$ be an arbitrary quiver.  The dimension vectors of
indecomposable representations of $Q$ correspond to positive roots
of the root system associated to the underlying graph of $Q$ (and
are thus independent of the orientation of the arrows of $Q$). The
correspondence is given by
\[
d_V \mapsto \sum_{i \in Q_0} d_V(i) \alpha_i.
\]
\end{theo}

Note that in Kac's Theorem, it is not asserted that the
isomorphism classes are in one-to-one correspondence with the
roots as in the finite case considered in Gabriel's theorem.  It
turns out that in the general case, dimension vectors for which
there is exactly one isomorphism class correspond to real roots
while imaginary roots correspond to dimension vectors for which
there are families of representations.

\begin{example}
Let $Q$ be the quiver of type $A_n$, oriented as follows.
\begin{center}
\psfrag{n2}{$n\! -\! 2$} \psfrag{n1}{$n\! -\! 1$} \psfrag{n}{$n$}
\psfrag{r1}{$\rho_1$} \psfrag{r2}{$\rho_2$}
\psfrag{rn2}{$\rho_{n-2}$} \psfrag{rn1}{$\rho_{n-1}$}
\includegraphics[width=0.4\textwidth]{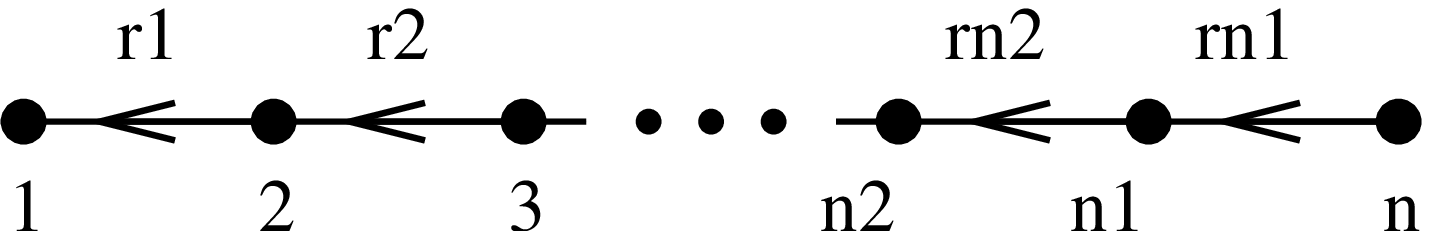}
\end{center}
It is known that the set of positive roots of the simple Lie
algebra of type $A_n$ is
\[
\left\{\left. \sum_{i=j}^l \alpha_i\ \right|\ 1 \le j \le l \le
n\right\} \sqcup \{0\}.
\]
The zero root corresponds to the trivial representation.  The root
$\sum_{i=j}^l \alpha_i$ for some $1 \le j \le l \le n$ corresponds
to the unique (up to isomorphism) representation $V$ with
\[
V_i = \begin{cases} k \text{ if } j \le i \le l \\
0 \text{ otherwise} \end{cases} \\
\]
and
\[
V_{\rho_i} = \begin{cases} 1 \text{ if } j \le i \le l-1 \\
0 \text{ otherwise } \end{cases}.
\]
\end{example}

\begin{example}
Let $Q$ be the quiver of type ${\widehat A}_n$, with all arrows
oriented in the same direction (for instance, counter-clockwise).
The positive root $\sum_{i=0}^n \alpha_i$ (where $\{0,1,2,\dots,
n\}$ are the vertices of the quiver) is imaginary.  There is a one
parameter family of isomorphism classes of indecomposable
representations where the maps assigned to each arrow are
non-zero.  The parameter is the composition of the maps around the
loop.
\end{example}

If a quiver $Q$ has no oriented cycles, then the only simple
$kQ$-modules are the modules $S^i$ for $i \in Q_0$ where
\[
S^i_j = \begin{cases} k \text{ if } i = j \\
0 \text{ if } i \ne j \end{cases}
\]
and $S^i_\rho = 0$ for all $\rho \in Q_1$.


\section{Ringel-Hall Algebras}

Let $k$ be the finite field $\mathbb{F}_q$ with $q$ elements and
let $Q$ be a quiver with no oriented cycles.  Let $\P$ be the set
of all isomorphism classes of $kQ$-modules which are finite as
sets (since $k$ is finite-dimensional, these are just the quiver
representations we considered above).  Let $A$ be a commutative
integral domain containing $\Z$ and elements $v, v^{-1}$ such that
$v^2 = q$.  The \emph{Ringel-Hall algebra} $H=H_{A,v}(kQ)$ is the
free $A$-module with basis $\{[V]\}$ indexed by the isomorphism
classes of representations of the quiver $Q$, with an $A$-bilinear
multiplication defined by
\[
[V^1] \cdot [V^2] = v^{\left<\dim V^1, \dim V^2\right>} \sum_V
g^V_{V^1,V^2} [V].
\]
Here $\left<\dim V^1, \dim V^2\right>$ is the Euler form and
$g^V_{V^1,V^2}$ is the number of submodules $W$ of $V$ such that
$V/W \cong V^1$ and $W \cong V^2$.  $H$ is an associative $\Z_{\ge
0}^{Q_0}$-graded algebra, with identity element $[0]$, the
isomorphism class of the trivial representation. The grading $H =
\oplus_{\alpha} H_\alpha$ is given by letting $H_\alpha$ be the
$A$-span of the set of isomorphism classes $[V]$ such that $\dim V
= \alpha$.

Let $C = C_{A,v}(kQ)$ be the $A$-subalgebra of $H$ generated by
the isomorphism classes $[S^i]$ of the simple $kQ$-modules.  $C$
is called the \emph{composition algebra}.  If the underlying graph
of $Q$ is of finite type, then $C = H$.

Now let $\mathcal{K}$ be a set of finite fields $k$ such that the
set $\{|k|\ |\ k \in \mathcal{K}\}$ is infinite. Let $A$ be an
integral domain containing $\Q$ and, for each $k \in \mathcal{K}$,
an element $v_k$ such that $v_k^2 = |k|$.  For each $k \in
\mathcal{K}$ we have the corresponding composition algebra $C_k$,
generated by the elements $[^kS^i]$ (here we make the field $k$
explicit). Now let $C$ be the subring of $\prod_{k \in
\mathcal{K}} C_k$ generated by $\Q$ and the elements
\begin{gather*}
t = (t_k)_{k \in \mathcal{K}},\ t_k = v_k, \\
t^{-1} = (t^{-1}_k)_{k \in \mathcal{K}},\ t^{-1}_k = (v^k)^{-1}, \\
u^i = (u^i_k)_{k \in \mathcal{K}},\ u^i_k = [^kS^i],\ i \in Q_0.
\end{gather*}
Now, $t$ lies in the center of $C$ and if $p(t)=0$ for some
polynomial $p$, then $p$ must be the zero polynomial since the set
of $v_k$ is infinite.  Thus we may think of $C$ as the $A$-algebra
generated by the $u^i$, $i \in Q_0$, with $A=\Q[t,t^{-1}]$ and $t$
an indeterminate.  Let $C^* = \Q(t) \otimes_A C$.  We call $C^*$
the \emph{generic composition algebra}.

Let $\g$ be the Kac-Moody algebra associated to the Cartan matrix
of the quiver $Q$ and let $U$ be the quantum group associated by
Drinfeld and Jimbo to $\g$.  It has a triangular decomposition $U
= U^- \otimes U^0 \otimes U^+$. Specifically, $U^+$ is the
$\Q(t)$-algebra with generators $E_i$, $i \in Q_0$ and relations
\[
\sum_{p=0}^{1-c_{ij}} (-1)^p \begin{bmatrix} 1-c_{ij} \\ p
\end{bmatrix} E_i^p E_j E_i^{1-c_{ij}-p},\ i \ne j,
\]
where $c_{ij}$ are the entries of the Cartan matrix and
\begin{gather*}
\begin{bmatrix} m \\ p \end{bmatrix} = \frac{[m]!}{[p]! [m-p]!},
\\
[n] = \frac{t^n - t^{-n}}{t-t^{-1}},\ [n]! = [1][2]\dots[n].
\end{gather*}

\begin{theo}
\label{theo:ringel-hall}
There is a $\Q(t)$-algebra isomorphism $C^* \to U^+$ sending $u_i
\mapsto E_i$ for all $i \in Q_0$.
\end{theo}
The proof of Theorem~\ref{theo:ringel-hall} is due to Ringel in
the case that the underlying graph of $Q$ is of finite or affine
type. The more general case presented here is due to Green.

All of the Kac-Moody algebras considered so far have been
simply-laced.  That is, their Cartan matrices are symmetric. There
is a way to deal with non-simply-laced Kac-Moody algebras using
\emph{species}.  We will not treat this subject in this article.


\section{Quiver Varieties}

One can use varieties associated to quivers to yield a geometric
realization of the upper half of the universal enveloping algebra
of a Kac-Moody algebra $\g$ and its irreducible highest-weight
representations.

\subsection{Lusztig's Quiver Varieties}

We first introduce the quiver varieties, first defined by Lusztig,
which yield a geometric realization of the upper half
$\mathcal{U}^+$ of the universal enveloping algebra of a
simply-laced Kac-Moody algebra $\g$.  Let $Q = (Q_0,Q_1)$ be the
quiver whose vertices $Q_0$ are the vertices of the Dynkin diagram
of $\g$ and whose set of arrows $Q_1$ consists of all the edges of
the Dynkin diagram with \emph{both} orientations. By definition,
$\mathcal{U}^+$ is the $\Q$-algebra defined by generators $e_i$, $i
\in Q_0$, subject to the Serre relations
\[
\sum_{p=0}^{1-c_{ij}} (-1)^p \begin{pmatrix} 1-c_{ij} \\ p
\end{pmatrix} e^p_i e_j e^{1-c_{ij}-p}
= 0
\]
for all $i \ne j$ in $Q_0$, where $c_{ij}$ are the entries of the
Cartan matrix associated to $Q$. For any $\nu = \sum_{i \in Q_0}
\nu_i i$, $\nu_i \in \mathbb{N}$, let $\mathcal{U}_\nu^+$ be the
subspace of $\mathcal{U}^+$ generated by the monomials $e_{i_1}
e_{i_2} \dots e_{i_n}$ for various sequences $i_1, i_2, \dots,
i_n$ in which $i$ appears $\nu_i$ times for each $i \in Q_0$. Thus
$\mathcal{U}^+ = \oplus_\nu \mathcal{U}_\nu^+$.  Let
$\mathcal{U}_\Z^+$ be the subring of $\mathcal{U}^+$ generated by
the elements $e_i^p/p!$ for $i \in Q_0$, $p \in \mathbb{N}$.  Then
$\mathcal{U}_\Z^+ = \oplus_\nu \mathcal{U}^+_{\Z,\nu}$ where
$\mathcal{U}^+_{\Z,\nu} = \mathcal{U}^+_\Z \cap
\mathcal{U}^+_\nu$.

We define the involution $\bar{\ }: Q_1 \to Q_1$ to be the
function which takes $\rho \in Q_1$ to the element of $Q_1$
consisting of the same edge with opposite orientation. An
\emph{orientation} of our graph/quiver is a choice of a subset
$\Omega \subset Q_1$ such that $\Omega \cup \bar{\Omega} = Q_1$
and $\Omega \cap \bar{\Omega} = \emptyset$.

Let $\mathcal{V}$ be the category of finite-dimensional
$Q_0$-graded vector spaces $\mathbf{V} = \oplus_{i
  \in Q_0} \mathbf{V}_i$ over $\C$ with morphisms being linear maps
respecting the grading.  Then $\mathbf{V} \in \mathcal{V}$ shall
denote that $\mathbf{V}$ is an object of $\mathcal{V}$.  The
dimension of $\mathbf{V} \in \mathcal{V}$ is given by $\mathbf{v}
= \dim \mathbf{V} = (\dim \mathbf{V}_0, \dots, \dim
\mathbf{V}_n)$.

Given $\mathbf{V} \in \mathcal{V}$, let $E_V$ be the space of
representations of $Q$ with underlying vector space $\mathbf{V}$.
That is,
\[
\mathbf{E_V} = \bigoplus_{\rho \in Q_1} \Hom
(\mathbf{V}_{t(\rho)}, \mathbf{V}_{h(\rho)}).
\]
For any subset $Q_1'$ of $Q_1$, let $\mathbf{E}_{\mathbf{V},
Q_1'}$ be the subspace of $\mathbf{E_V}$ consisting of all vectors
$x = (x_\rho)$ such that $x_\rho=0$ whenever $\rho \not\in Q_1'$.
The algebraic group $G_\mathbf{V} = \prod_i \Aut(\mathbf{V}_i)$
acts on $\mathbf{E_V}$ and $\mathbf{E}_{\mathbf{V}, Q_1'}$ by
\[
(g,x) = ((g_i), (x_\rho)) \mapsto gx = (x'_\rho) = (g_{h(\rho)}
x_\rho g_{t(\rho)}^{-1}).
\]

Define the function $\varepsilon : Q_1 \to \{-1,1\}$ by
$\varepsilon (\rho) = 1$ for all $\rho \in \Omega$ and
$\varepsilon(\rho) = -1$ for all $\rho \in {\bar{\Omega}}$.
Let
$\left<\cdot,\cdot\right>$ be the nondegenerate,
$G_\mathbf{V}$-invariant, symplectic form on $\mathbf{E_V}$ with
values in $\C$ defined by
\[
\left<x,y\right> = \sum_{\rho \in Q_1} \varepsilon(\rho) \tr
(x_\rho y_{\bar{\rho}}).
\]
Note that $\mathbf{E_V}$ can be considered as the cotangent space
of $\mathbf{E}_{\mathbf{V}, \Omega}$ under this form.

The moment map associated to the $G_{\mathbf{V}}$-action on the
symplectic vector space $\mathbf{E_V}$ is the map $\psi :
\mathbf{E_V} \to \mathbf{gl_V} = \prod_i \End \mathbf{V}_i$, the
Lie algebra of $GL_\mathbf{V}$, with $i$-component $\psi_i :
\mathbf{E_V} \to \End \mathbf{V}_i$ given by
\[
\psi_i(x) = \sum_{\rho \in Q_1,\, h(\rho)=i} \varepsilon(\rho)
x_\rho x_{\bar{\rho}} .
\]

\begin{defin}
\label{def:nilpotent}
An element $x \in \mathbf{E_V}$ is said to be \emph{nilpotent} if
there exists an $N \ge 1$ such that for any sequence $\rho_1,
\rho_2, \dots, \rho_N$ in $H$ satisfying $t (\rho_1) = h
(\rho_2)$, $t (\rho_2) = h (\rho_3)$, \dots, $t (\rho_{N-1}) = h
(\rho_N)$, the composition $x_{\rho_1} x_{\rho_2} \dots x_{\rho_N}
: \mathbf{V}_{t (\rho_N)} \to
  \mathbf{V}_{h (\rho_1)}$ is zero.
\end{defin}

\begin{defin} Let $\Lambda_\mathbf{V}$ be the set of all
  nilpotent elements $x \in \mathbf{E_V}$ such that $\psi_i(x) = 0$
  for all $i \in I$.
\end{defin}

A subset of an algebraic variety is said to be
\emph{constructible} if it is obtained from subvarieties from a
finite number of the usual set-theoretic operations.  A function
$f: A \to \Q$ on an algebraic variety $A$ is said to be a
\emph{constructible function} if $f^{-1}(a)$ is a constructible
set for all $a \in \Q$ and is empty for all but finitely many $a$.
Let $M(\Lambda_\mathbf{V})$ denote the $\Q$-vector space of all
constructible functions on $\Lambda_\V$.  Let ${\widetilde
M}(\Lambda_\V)$ denote the $\Q$-subspace of $M(\Lambda_\V)$
consisting of those functions that are constant on any
$G_\V$-orbit in $\Lambda_\V$.

Let $\V, \V', \V'' \in \mathcal{V}$ such that $\dim \V = \dim \V'
+ \dim \V''$.  Now, suppose that $\mathbf{S}$ is an $I$-graded
subspace of $\mathbf{V}$. For $x \in \Lambda_\mathbf{V}$ we say
that $\mathbf{S}$ is \emph{$x$-stable} if $x(\mathbf{S}) \subset
\mathbf{S}$.  Let $\Lambda_{\V;\V',\V''}$ be the variety
consisting of all pairs $(x,\S)$ where $x \in \Lambda_\V$ and $\S$
is an $I$-graded $x$-stable subspace of $\V$ such that $\dim \S =
\dim \V''$.  Now, if we fix some isomorphisms $\V/\S \cong \V'$,
$\S \cong \V''$, then $x$ induces elements $x' \in \Lambda_{\V'}$
and $x'' \in \Lambda_{\V''}$.  We then have the maps
\[
\Lambda_{\V'} \times \Lambda_{\V''} \stackrel{p_1}{\longleftarrow}
\Lambda_{\V;\V',\V''} \stackrel{p_2}{\longrightarrow} \Lambda_{\V}
\]
where $p_1(x,\S) = (x',x'')$, $p_2(x,\S) = x$.

For a holomorphic map $\pi$ between complex varieties $A$ and $B$,
let $\pi_!$ denote the map between the spaces of constructible
functions on $A$ and $B$ given by
\[
(\pi_! f)(y) = \sum_{a \in \Q} a \chi (\pi^{-1}(y) \cap f^{-1}(a)).
\]
Let $\pi^*$ be the pullback map from functions on $B$ to functions
on $A$ acting as $\pi^* f(y) = f(\pi(y))$. We then define a map
\begin{equation}
\label{eq:func-mult}
{\widetilde M}(\Lambda_{\V'}) \times {\widetilde
M}(\Lambda_{\V''}) \to {\widetilde M}(\Lambda_{\V})
\end{equation}
by $(f',f'') \mapsto f' \ast f''$ where
\[
f' \ast f'' = (p_2)_! p_1^* (f' \times f'').
\]
Here $f' \times f'' \in {\widetilde M}(\Lambda_{\V'} \times
\Lambda_{\V''})$ is defined by $(f' \times f'')(x',x'') = f'(x')
f''(x'')$.  The map \eqref{eq:func-mult} is bilinear and defines
an associatve $\Q$-algebra structure on $\oplus_\nu {\widetilde
M}(\Lambda_{\V^{\nu}})$ where $\V^\nu$ is the object of
$\mathcal{V}$ defined by $\V^\nu_i = \C^{\nu_i}$.

There is a unique algebra homomorphism $\kappa : \mathcal{U}^+ \to
\oplus_\nu {\widetilde M}(\Lambda_{\V^\nu})$ such that
$\kappa(e_i)$ is the function on the point $\Lambda_{\V^i}$ with
value 1.  Then $\kappa$ restricts to a map $\kappa_{\nu} :
\mathcal{U}_\nu^+ \to {\widetilde M}(\Lambda_{\V^\nu})$.  It can
be shown that $\kappa_{pi}(e_i^p/p!)$ is the function 1 on the
point $\Lambda_{\V^{pi}}$ for $i \in Q_0$, $p \in \Z_{\ge 0}$.

Let ${\widetilde M}_\Z(\Lambda_\V)$ be the set of all functions in
${\widetilde M}(\Lambda_\V)$ that take on only integer values. One
can show that if $f' \in {\widetilde M}_\Z(\Lambda_{\V'})$ and
$f'' \in {\widetilde M}_\Z(\Lambda_{\V''})$ then $f' \ast f'' \in
{\widetilde M}_\Z(\Lambda_\V)$ in the setup of
\eqref{eq:func-mult}.  Thus $\kappa_\nu(\mathcal{U}^+_{\Z, \nu})
\subseteq {\widetilde M}_\Z(\Lambda_{\V^\nu})$.

Let $\Irr \Lambda_\V$ denote the set of irreducible components of
$\Lambda_\V$.  The following proposition was conjectured by
Lusztig and proved by him in the affine (and finite) case.  The
general case was proved by Kashiwara and Saito.

\begin{prop}
For any $\nu \in (\Z_{\ge 0})^{Q_0}$, we have $\dim
\mathcal{U}_\nu^+ = \# \Irr \Lambda_{\V^\nu}$.
\end{prop}

We then have the following important result due to Lusztig.
\begin{theo}
Let $\nu \in (\Z_{\ge 0})^{Q_0}$.  Then
\begin{enumerate}
\item For any $Z \in \Irr \Lambda_{\V^\nu}$, there exists a unique
$f_Z \in \kappa_\nu(\mathcal{U}^+_{\Z,\nu})$ such that $f_Z$ is
equal to one on an open dense subset of $Z$ and equal to zero on an
open dense subset of $Z' \in \Irr \Lambda_{\V^\nu}$ for all $Z' \ne
Z$.

\item $\{f_Z\ |\ Z \in \Irr \Lambda_{\V^\nu}\}$ is a $\Q$-basis of
$\kappa_\nu(\mathcal{U}^+_\nu)$.

\item $\kappa_\nu : \mathcal{U}^+_\nu \to
\kappa_\nu(\mathcal{U}^+_\nu)$ is an isomorphism.

\item Define $[Z] \in \mathcal{U}^+_\nu$ by $\kappa_\nu([Z]) =
f_Z$.  Then $B_\nu = \{[Z]\ |\ Z \in \Irr \Lambda_{\V^\nu}\}$ is a
$\Q$-basis of $\mathcal{U}^+_\nu$.

\item $\kappa_\nu(\mathcal{U}^+_{\Z,\nu}) =
\kappa_\nu(\mathcal{U}^+_\nu) \cap {\widetilde
M}_\Z(\Lambda_{\V^\nu})$.

\item $B_\nu$ is a $\Z$-basis of $\mathcal{U}^+_{\Z,\nu}$.
\end{enumerate}
\end{theo}

From this theorem, we see that $B = \sqcup_\nu B_\nu$ is a
$\Q$-basis of $\mathcal{U}^+$, which is called the
\emph{semicanonical basis}.  This basis has many remarkable
properties.  One of these properties is as follows.  Via the
algebra involution of the entire universal enveloping algebra
$\mathcal{U}$ of $\g$ given on the Chevalley generators by $e_i
\mapsto f_i$, $f_i \mapsto e_i$ and $h \mapsto -h$ for $h$ in the
Cartan subalgebra of $\g$, one obtains from the results of this
section a semicanonical basis of $\mathcal{U}^-$, the lower half
of the universal enveloping algebra of $\g$.  For any irreducible
highest weight integrable representation $V$ of $\mathcal{U}$ (or,
equivalently, $\g$), let $v \in V$ be a non-zero highest weight
vector. Then the set
\[
\{bv\ |\ b \in B,\, bv \ne 0\}
\]
is a $\Q$-basis of $V$, called the \emph{semicanonical basis} of
$V$.  Thus the semicanonical basis of $\mathcal{U}^-$ is
simultaneously compatible will all irreducible highest weight
integrable modules.  There is also a way to define the
semicanonical basis of a representation directly in a geometric
way. This is the subject of the next subsection.

One can also obtain a geometric realization of the upper part
$U^+$ of the quantum group in a similar manner using perverse
sheaves instead of constructible functions.  This construction
yields the \emph{canonical basis} which also has many remarkable
properties and is closely related to the theory of crystal bases.


\subsection{Nakajima's Quiver Varieties}

We introduce here a description of the quiver varieties first
presented by Nakajima.  They yield a geometric realization of the
irreducible highest weight representations of a simply-laced
Kac-Moody algebra.

\begin{defin}
\label{def:lambda}
For $\mathbf{v}, \mathbf{w} \in \Z_{\ge 0}^I$, choose $I$-graded
vector spaces $\mathbf{V}$ and $\mathbf{W}$ of graded dimensions
$\mathbf{v}$ and $\mathbf{w}$ respectively.  Then define
\[
\Lambda \equiv \Lambda(\mathbf{v},\mathbf{w}) = \Lambda_\mathbf{V}
\times \sum_{i \in I} \Hom (\mathbf{V}_i, \mathbf{W}_i).
\]
\end{defin}

\begin{defin}
\label{def:lambda-stable}
Let $\Lambda^{\text{st}} =
\Lambda(\mathbf{v},\mathbf{w})^{\text{st}}$ be the set of all $(x,
t) \in \Lambda(\mathbf{v},\mathbf{w})$ satisfying the following
condition:  If $\mathbf{S}=(\mathbf{S}_i)$ with $\mathbf{S}_i
\subset \mathbf{V}_i$ is $x$-stable and $t_i(\mathbf{S}_i) = 0$
for $i \in I$, then $\mathbf{S}_i = 0$ for $i \in I$.
\end{defin}

The group $G_\mathbf{V}$ acts on $\Lambda(\mathbf{v},\mathbf{w})$
via
\[
(g,(x,t)) \mapsto ((g_{h (\rho)} x_\rho g_{t (\rho)}^{-1}), (t_i
g_i^{-1})).
\]
and the stabilizer of any point of
$\Lambda(\mathbf{v},\mathbf{w})^{\text{st}}$ in $G_{\mathbf{V}}$
is trivial.  We then make the following definition.
\begin{defin}
\label{def:L}
Let $\mathcal{L} \equiv \mathcal{L}(\mathbf{v},\mathbf{w}) =
\Lambda(\mathbf{v},\mathbf{w})^{\text{st}} / G_{\mathbf{V}}$.
\end{defin}

Let $\mathbf{w, v, v', v''} \in \Z_{\ge 0}^I$ be such that
$\mathbf{v} = \mathbf{v'} + \mathbf{v''}$.  Consider the maps
\begin{equation}
\label{eq:diag_action}
\Lambda(\mathbf{v}'',\mathbf{0}) \times
\Lambda(\mathbf{v}',\mathbf{w}) \stackrel{p_1}{\leftarrow}
\mathbf{\tilde F (v,w;v'')} \stackrel{p_2}{\rightarrow}
\mathbf{F(v,w;v'')} \stackrel{p_3}{\rightarrow}
\Lambda(\mathbf{v},\mathbf{w}),
\end{equation}
where the notation is as follows.  A point of
$\mathbf{F(v,w;v'')}$ is a point $(x,t) \in
\Lambda(\mathbf{v},\mathbf{w})$ together with an $I$-graded,
$x$-stable subspace $\mathbf{S}$ of $\mathbf{V}$ such that $\dim
\mathbf{S} = \mathbf{v'} = \mathbf{v} - \mathbf{v''}$.  A point of
$\mathbf{\tilde
  F (v,w;v'')}$ is a point $(x,t,\mathbf{S})$ of $\mathbf{F(v,w;v'')}$
together with a collection of isomorphisms $R'_i : \mathbf{V}'_i
\cong \mathbf{S}_i$ and $R''_i : \mathbf{V}''_i \cong \mathbf{V}_i
/ \mathbf{S}_i$ for each $i \in I$.  Then we define
$p_2(x,t,\mathbf{S}, R',R'') = (x,t,\mathbf{S})$,
$p_3(x,t,\mathbf{S}) = (x,t)$ and $p_1(x,t,\mathbf{S},R',R'') =
(x'',x',t')$ where $x'', x', t'$ are determined by
\begin{align*}
R'_{h(\rho)} x'_\rho &= x_\rho R'_{t(\rho)} :
\mathbf{V}'_{t(\rho)} \to
\mathbf{S}_{h(\rho)}, \\
t'_i &= t_i R'_i : \mathbf{V}'_i \to \mathbf{W}_i \\
R''_{h(\rho)} x''_\rho &= x_\rho R''_{t(\rho)} :
\mathbf{V}''_{t(\rho)} \to \mathbf{V}_{h(\rho)} /
\mathbf{S}_{h(\rho)}.
\end{align*}
It follows that $x'$ and $x''$ are nilpotent.

\begin{lem}
One has
\[
(p_3 \circ p_2)^{-1} (\Lambda(\mathbf{v},\mathbf{w})^{\text{st}})
\subset p_1^{-1} (\Lambda(\mathbf{v}'',\mathbf{0}) \times
\Lambda(\mathbf{v}',\mathbf{w})^{\text{st}}).
\]
\end{lem}

Thus, we can restrict \eqref{eq:diag_action} to
$\Lambda^{\text{st}}$, forget the
$\Lambda(\mathbf{v}'',\mathbf{0})$-factor and consider the
quotient by $G_\mathbf{V}$ and $G_\mathbf{V'}$.  This yields the
diagram
\begin{equation}
\label{eq:diag_action_mod}
\mathcal{L}(\mathbf{v'}, \mathbf{w}) \stackrel{\pi_1}{\leftarrow}
\mathcal{F}(\mathbf{v}, \mathbf{w}; \mathbf{v} - \mathbf{v'})
\stackrel{\pi_2}{\rightarrow} \mathcal{L}(\mathbf{v}, \mathbf{w}),
\end{equation}
where
\begin{gather*}
\mathcal{F}(\mathbf{v}, \mathbf{w}, \mathbf{v} - \mathbf{v'}) \\
\stackrel{\text{def}}{=} \{ (x,t,\mathbf{S}) \in
\mathbf{F(v,w;v-v')}\,
  |\, (x,t) \in \Lambda(\mathbf{v},\mathbf{w})^{\text{st}} \} / G_\mathbf{V}.
\end{gather*}

Let $M(\mathcal{L}(\mathbf{v}, \mathbf{w}))$ be the vector space
of all constructible functions on $\mathcal{L}(\mathbf{v},
\mathbf{w})$. Then define maps
\begin{gather*}
h_i : M(\mathcal{L}(\mathbf{v}, \mathbf{w})) \to
M(\mathcal{L}(\mathbf{v}, \mathbf{w})) \\
e_i : M(\mathcal{L}(\mathbf{v}, \mathbf{w})) \to
M(\mathcal{L}(\mathbf{v} - \mathbf{e}^i, \mathbf{w}))  \\
f_i : M(\mathcal{L}(\mathbf{v} - \mathbf{e}^i, \mathbf{w})) \to
M(\mathcal{L}(\mathbf{v}, \mathbf{w})).
\end{gather*}
by
\begin{gather*}
h_i f = u_i f, \\
e_i f = (\pi_1)_! (\pi_2^* f), \\
f_i g = (\pi_2)_! (\pi_1^* g).
\end{gather*}
Here
\[
\mathbf{u} = {^t(u_0, \dots, u_n)} = \mathbf{w} - C \mathbf{v}
\]
where $C$ is the Cartan matrix of $\mathfrak{g}$ and we are using
diagram~\eqref{eq:diag_action_mod} with $\mathbf{v}' = \mathbf{v} -
\mathbf{e}^i$ where $\mathbf{e}^i$ is the vector whose components
are given by $\mathbf{e}^i_j = \delta_{ij}$.

Now let $\varphi$ be the constant function on
$\mathcal{L}(\mathbf{0}, \mathbf{w})$ with value 1.  Let
$L(\mathbf{w})$ be the vector space of functions generated by
acting on $\varphi$ with all possible combinations of the
operators $f_i$.  Then let $L(\mathbf{v},\mathbf{w}) =
M(\mathcal{L}(\mathbf{v}, \mathbf{w})) \cap L(\mathbf{w})$.

\begin{prop}
The operators $e_i$, $f_i$, $h_i$ on $L(\mathbf{w})$ provide it
with the structure of the irreducible highest weight integrable
representation of $\mathfrak{g}$ with highest weight $\sum_{i \in
Q_0} \mathbf{w}_i \omega_i$. Each summand of the decomposition
$L(\mathbf{w}) = \bigoplus_\mathbf{v} L(\mathbf{v}, \mathbf{w})$
is a weight space with weight $\sum_{i \in Q_0} \mathbf{w}_i
\omega_i - \mathbf{v}_i \alpha_i$.  Here the $\omega_i$ and
$\alpha_i$ are the fundamental weights and simple roots of
$\mathfrak{g}$ respectively.
\end{prop}

Let $Z \in \Irr \mathcal{L}(\mathbf{v}, \mathbf{w})$ and define a
linear map $T_Z : L(\mathbf{v}, \mathbf{w}) \to \C$ that
associates to a constructible function $f \in L(\mathbf{v},
\mathbf{w})$ the (constant) value of $f$ on a suitable open dense
subset of $Z$. The fact that $L(\mathbf{v}, \mathbf{w})$ is
finite-dimensional allows us to take such an open set on which
\emph{any} $f \in L(\mathbf{v}, \mathbf{w})$ is constant.  So we
have a linear map
\[
\Phi : L(\mathbf{v}, \mathbf{w}) \to \C^{\Irr
\mathcal{L}(\mathbf{v},
  \mathbf{w})}.
\]
Then we have the following proposition.

\begin{prop}
\label{prop:func_irrcomp_isom}
The map $\Phi$ is an isomorphism; for any $Z \in \Irr
\mathcal{L}(\mathbf{v}, \mathbf{w})$, there is a unique function
$g_Z \in L(\mathbf{v}, \mathbf{w})$ such that for some open dense
subset $O$ of $Z$ we have $g_Z|_O = 1$ and for some closed
$G_\mathbf{V}$-invariant subset $K \subset \mathcal{L}(\mathbf{v},
\mathbf{w})$ of dimension $< \dim \mathcal{L}(\mathbf{v},
\mathbf{w})$ we have $g_Z=0$ outside $Z \cup K$.  The functions
$g_Z$ for $Z \in \Irr \Lambda(\mathbf{v},\mathbf{w})$ form a basis
of $L(\mathbf{v},\mathbf{w})$.
\end{prop}


\section{Further Reading}

We have given here a brief overview of some topics related to
finite-dimensional algebras and quivers.  There is much more to be
found in the literature.  For basics on associative algebras and
their representations, the reader may consult introductory texts on
abstract algebra such as \cite{Lang}.  For further results (and
their proofs) on Ringel-Hall algebras see the papers
\cite{Rin90b,Rin90a,Rin91,Rin93,Rin94} of Ringel and Green's paper
\cite{Gre} and the references cited therein.  The reader interested
in species, which extend many of these results to non-simply-laced
Lie algebras, should consult \cite{DR76}.

The book \cite{L93} covers the quiver varieties of Lusztig and
canonical bases.  Canonical bases are closely related to crystal
bases and crystal graphs (see \cite{HK} for an overview of these
topics).  In fact, the set of irreducible components of the quiver
varieties of Lusztig and Nakajima can be endowed with the
structure of a crystal graph in a purely geometric way (see
\cite{KS97,S02}).  Many results on Nakajima's quiver varieties can
be found in the original papers \cite{N94,N98}.  The overview
article \cite{N96} is also useful.

Quiver varieties can also be used to give geometric realizations
of tensor products of representations (see
\cite{Mal02,Mal03,Nak01a,Sav03a}) and finite-dimensional
representations of quantum affine Lie algebras (see
\cite{Nak01b}).  This is just a select few of the many
applications of quiver varieties.  Much more can be found in the
literature.


\bibliographystyle{abbrv}
\bibliography{biblist}

\end{multicols}

\end{document}